\documentclass[preprint,10pt]{elsarticle}

\usepackage{amssymb}
\usepackage{amsthm}
\usepackage{amsmath}
 \newtheorem{thm}{Theorem}[section]
 \newtheorem{cor}[thm]{Corollary}
 \newtheorem{lem}[thm]{Lemma}
 
 \theoremstyle{definition}
 
 \theoremstyle{remark}
 \newtheorem{rem}[thm]{Remark}

\journal{arxiv}

\begin{document}

\begin{frontmatter}



\title{Some $L^p$ rigidity results for complete  manifolds with harmonic curvature}


\author[HPF]{Hai-Ping Fu}
\ead{mathfu@126.com}

\author[LQX]{Li-Qun Xiao}
\ead{xiaoliqun@ncu.edu.cn}


\fntext[fn1]{Supported  by National Natural Science Foundations of China (11261038, 11361041),  Jiangxi Province
Natural Science Foundation of China (20132BAB201005).}

\address[HPF]{Department of Mathematics,  Nanchang University, Nanchang 330031, P.
R. China}

    \address[LQX]{Department of management science and Engineering,  Nanchang University, Nanchang 330031, P.
R. China}

\begin{abstract}
Let $(M^n, g)(n\geq3)$ be an $n$-dimensional complete  Riemannian manifold with harmonic curvature and positive Yamabe constant. Denote by  $R$ and $\mathring{Rm}$ the scalar curvature and the trace-free Riemannian curvature tensor  of $M$, respectively. The main result of this paper states that $\mathring{Rm}$ goes to
zero uniformly at infinity if for $p\geq \frac n2$, the $L^{p}$-norm of  $\mathring{Rm}$ is finite. Moreover, If $R$ is positive, then $(M^n, g)$ is compact. As applications,  we prove that $(M^n, g)$ is isometric to a spherical
space form if for $p\geq \frac n2$, $R$ is positive and the $L^{p}$-norm of $\mathring{Rm}$ is pinched in $[0,C_1)$, where $C_1$ is
an explicit positive constant depending only on $n, p$, $R$ and the Yamabe constant.

In particular, we prove an $L^{p}(\frac n2\leq p<\frac{n-2}{2}(1+\sqrt{1-\frac4n}))$-norm of $\mathring{Ric}$  pinching theorem
for   complete, simply connected,  locally conformally flat Riemannian $n(n\geq 6)$-manifolds with constant
negative scalar curvature.

 We give an isolation theorem of the trace-free Ricci curvature tensor of  compact  locally conformally flat
Riemannian $n$-manifolds with constant positive scalar curvature, which improves Thereom 1.1 and Corollary 1 of E. Hebey  and  M. Vaugon \cite{{HV}}. This rsult is sharped, and we can precisely characterize the case of equality.
\end{abstract}


\end{frontmatter}


\section{Introduction and main results}
\label{}

Recall that an $n$-dimensional Riemannian manifold $(M^n, g)$  is said to be  a manifold with harmonic curvature
if the divergence of its Riemannian curvature tensor $Rm$ vanishes, i.e., $\delta Rm=0$. In view of the second Bianchi identity, we know that $M$ has harmonic
curvature if and only if the Ricci tensor of $M$ is a Codazzi tensor. When $n\geq 3$, by the Bianchi identity, the scalar curvature is constant. Thus, every Riemannian manifold with parallel Ricci tensor has harmonic
curvature. Moreover, the constant curvature spaces, Einstein manifolds and the locally conformally
flat manifolds with constant scalar curvature are
also important examples of manifolds with harmonic curvature, however, the
converse does not hold (see \cite{{B}}, for example). According to the decomposition of the Riemannian curvature tensor, the metric with harmonic curvature is a natural candidate for this study  since one of the important problems in Riemannian geometry is to understand classes
of metrics that are, in some sense, close to being Einstein or having constant curvature. The another reason for this study on the metric with harmonic curvature is the fact that a Riemannian manifold has harmonic curvature if and
only if the Riemannian connection is a solution of the Yang-Mills equations on
the tangent bundle \cite{Bo}. In recent years, the complete manifolds  with
harmonic curvature have been studied in literature (e.g., \cite{{C},{Ch},{FDS},{FP},{HV},{IS},{K},{PRS},{S},{XZ},{TV}}).
In particular, G. Tian and J. Viaclovsky \cite{TV},  and   X. Chen and B. Weber \cite{CW} have obtained $\epsilon$-rigidity results for critical metric which relies on a Sobolev inequality and a integral bounds on the
curvature in
dimension $4$ and in higher dimension, respectively.
The curvature pinching phenomenon plays an important role in global differential
geometry. We are interested in  $L^p$
pinching problems for complete Riemannian manifold with harmonic curvature.

We  now introduce the definition of the Yamabe constant. Given a complete Riemannian n-manifold $M$, the Yamabe constant $Q(M)$ is defined by
$$Q(M)=\inf_{0\neq u\in C_0^{\infty}(M)}\frac{\int_{M}\left(|\nabla u|^2+\frac{(n-2)}{4(n-1)}R
u^{2} \right)}{\left(\int_{M} |u|^{\frac{2n}{n-2}}\right)^{\frac{n-2}{n}}},$$
where $R$ is the scalar curvature of $M$. The important works of Schoen, Trudinger and Yamabe showed that the infimum in the above on compact manifolds is always achieved (see \cite{{A}, {SY}}).
There are  complete noncompact Riemannian manifolds of negative scalar
curvature with positive Yamabe constant. For example, any simply connected
complete locally conformally flat manifold has positive Yamabe constant \cite{SY1}, and $Q(M)$ is always positive if $R$  vanishes \cite{D}. In
contrast with the noncompact case, the Yamabe constant of a given compact manifold
is determined by the sign of scalar curvature \cite{A}.

 Throughout this paper, we always assume
that $M$ is an $n$-dimensional complete Riemannian manifold with $n\geq3$. In this note, we obtain the following rigidity theorems.
\begin{thm}
Let $M$ be a complete noncompact Riemannian $n$-manifold
with harmonic curvature. Assume that $M$ has the positive Yamabe constant or satisfies the Sobolev inequality
$$\left(\int_{M}|f|^{\frac{2n}{n-2}}\right)^{\frac{n-2}{n}}\leq C_{S}\int_{M}|\nabla f|^2, \forall f\in C^{\infty}_{0}(M). $$
For $p\geq \frac n2$, if $\int_{M}|\mathring{Rm}|^{p}< +\infty,$ then, given any $\epsilon>0$ and any $x_0\in M$ there exists a
geodesic ball $B_{r}(x_0)$ with center $x_0$ and radius $r$ such that $|\mathring{Rm}|(x)<\epsilon$ for all
$x\in M\setminus B_{r}(x_0)$.
\end{thm}
\begin{thm}
Let $M$ be a complete  Riemannian $n$-manifold
with harmonic curvature and positive scalar curvature. Assume that $M$ has the positive Yamabe constant.
For $p\geq \frac n2$, if $\int_{M}|\mathring{Rm}|^{p}< +\infty,$ then $M$ must be compact.
\end{thm}
\begin{cor}
Let $M$ be a complete noncompact  Riemannian $n$-manifold
with harmonic curvature and nonnegative scalar curvature. Assume that $M$ has the positive Yamabe constant.
For $p\geq \frac n2$, if $\int_{M}|\mathring{Rm}|^{p}< +\infty,$ then $M$ must be scalar flat.
\end{cor}

\begin{thm}
Let $M$ be a complete  Riemannian $n$-manifold
with harmonic curvature and positive scalar curvature $R$. Assume that $M$ has the positive Yamabe constant.
For $p\geq \frac n2,$   if
$$\left(\int_{M}|\mathring{Rm}|^{p}\right)^{\frac 1p}<C_1, $$
where  \begin{equation*}C_1=
\begin{cases}\frac {Q(M)}{C(n)}, \quad 3\leq n\leq5 \ \text{and}\  p=\frac n2\\ \frac{2(6-n)p R}{4(n-1)(2p-n)C(n)}[\frac{4(n-1)(2p-n)Q(M)}{(6-n)nR}]^{\frac{n}{2p}}, \quad 3\leq n\leq5 \ \text{and}\  \frac n2<p<\frac {2n}{n-2}\\ \frac{R}{(n-1)C(n)}\left[\frac{4(n-1)Q(M)}{(n-2)R}\right]^{\frac {n}{2p}}, \quad 3\leq n\leq5 \ \text{and}\  p\geq\frac {2n}{n-2}, \text{and}\  n\geq6,
\end{cases}
\end{equation*}
and $C(n)$ is defined in Lemma 2.1,
then $M$ is isometric to a spherical
space form.
\end{thm}
\begin{rem} Some $L^{\frac{n}{2}}$ and $L^{n}$
trace-free Riemannian curvature pinching theorems have been shown by Kim \cite{K}, Chu \cite{C}, and Fu etc.\cite{FDS}, in which the constant $C$ is not explicit, respectively. When $p=\frac n2$  or $p=n$, the constant $C_1$ in Theorem 1.4 satisfies
\begin{equation*}C_1=
\begin{cases}\frac {Q(M)}{C(n)}, \quad 3\leq n\leq5 \\ \frac{4Q(M)}{(n-2)C(n)}, \quad n\geq6
\end{cases}\text{or}\   C_1=
\begin{cases}\frac {\sqrt{3Q(M)R}}{\sqrt{2}C(3)}, \quad  n=3 \\(\frac{4Q(M)}{(n-1)(n-2)})^{\frac 12}\frac {1}{C(n)}, \quad n\geq4.
\end{cases}
\end{equation*}
 For compact Einstien manifolds or compact conformally fiat manifolds, there are some corresponding results in \cite{{Ca2}, {FX},{HV},{IS},{S}}.
\end{rem}

\begin{thm}
Let $M^n(n\geq 10)$ be a complete  Riemannian $n$-manifold
with harmonic curvature and negative scalar scalar curvature $R$. Assume that $M$ has the positive Yamabe constant.
If $$\int_{M}|\mathring{Rm}|^{\gamma}<\infty, \text{for some}\quad \gamma\in (0, \frac{n(n-2)+\sqrt{n(n-2)(n^2-10n+8)}}{4(n-1)}).$$
Then there exists a small number $C$ such that if for $p\geq \frac n2$,
$$\int_{M}|\mathring{Rm}|^{p}<C,$$
then $M$ is a hyperbolic space form.
\end{thm}
\begin{rem} Theorem 1.6 can be considered as generalization of  Theorem 1.4 in \cite{FDS} and some result in \cite{K}.
\end{rem}
When $M$ is a  complete, simply connected, locally
conformally flat Riemannian $n$-manifold, $M$ satisfies the Sobolev inequality (see Corollary 3.2 in  \cite{H}). Based on $(9)$, using the same argument as in the proofs of Theorems 1.4 and 1.6, we generalize the result due to \cite{XZ} and \cite{FP}.
\begin{thm}
Let $M$ be a complete, simply connected,  locally conformally flat
Riemannian $n$-manifold with constant positive scalar curvature. For $p\geq \frac n2$,
if
$$\left(\int_{M}|\mathring{Ric}|^{p}\right)^{\frac 1p}<C_2,$$
where \begin{equation*}C_2=
\begin{cases}\frac{3\sqrt{6}}{4}\omega_3^{\frac 23}, \quad n=3 \ \text{and}\  p=\frac 32\\ [\frac{6(2p-3)}{R}]^{\frac{3}{2p}}\frac{\sqrt{6}pR}{12(2p-3)}\omega^{\frac 1p}_3, \quad n=3 \ \text{and}\  \frac 32< p<2\\   [\frac{n(n-1)}{R}]^{\frac{n}{2p}}\frac{R}{\sqrt{n(n-1)}}\omega^{\frac 1p}_n, \quad n=3 \ \text{and}\  p\geq2, \text{and}\  n\geq4,
\end{cases}
\end{equation*}
then $M$ is isometric to a sphere.
\end{thm}
\begin{rem} Theorem 1.8  improves the $L^{\frac{n}{2}}(n\geq6)$  trace-free Ricci curvature pinching theorem given by \cite{XZ} in dimension. The pinching constant in Theorem 1.8 is better  than the one in the $L^n$  trace-free Ricci curvature pinching theorem given by \cite{XZ}.
\end{rem}
\begin{thm}
Let $M^n(n\geq 5)$ be a complete, simply connected,  locally conformally flat
Riemannian $n$-manifold with constant negative scalar curvature. Assume that
$$\int_{M}|\mathring{Ric}|^{\gamma}<\infty, \text{for some}\quad \gamma\in (0, \frac{n-2}{2}(1+\sqrt{1-\frac4n})).$$
Then there exists a small number $C$ such that if for $p\geq \frac n2$,
$$\int_{M}|\mathring{Ric}|^{p}< C,$$
then $M$ is a hyperbolic space.
\end{thm}
\begin{cor}
Let $M^n(n\geq 6)$ be a complete, simply connected,  locally conformally flat
Riemannian $n$-manifold with constant negative   scalar curvature.
For some  $p\in [\frac n2, \frac{n-2}{2}(1+\sqrt{1-\frac4n}))$, there exists a small number $C$ such that if
$$\left(\int_{M}|\mathring{Ric}|^{p}\right)^{\frac 1p}< C,$$
then $M$ is a hyperbolic space. In particular, when $p=\frac{n}{2}$, there exists an explicit positive constant $C=\sqrt{n(n-1)}\omega_{n}^{\frac2n}.$
\end{cor}
\begin{rem} Theorems 1.8 and 1.9 improve the corresponding one in \cite{{FP},{XZ}}. Corollary 1.11  has been proved in \cite{FP}.
\end{rem}

Using the same argument as in the proof of Theorem 1.8, we obtain
\begin{thm}
Let $M$ be a compact  locally conformally flat
Riemannian $n$-manifold with constant positive scalar curvature. For $p\geq \frac n2$,
if
$$\left(\int_{M}|\mathring{Ric}|^{p}\right)^{\frac 1p}<C_3,$$
where \begin{equation*}C_3=
\begin{cases}\sqrt{6}Q(M), \quad n=3 \ \text{and}\  p=\frac 32\\ [\frac{8(2p-3)Q(M)}{R}]^{\frac{3}{2p}}\frac{\sqrt{6}pR}{12(2p-3)}, \quad n=3 \ \text{and}\  \frac 32< p<2\\   [\frac{4(n-1)Q(M)}{(n-2)R}]^{\frac{n}{2p}}\frac{R}{\sqrt{n(n-1)}}, \quad n=3 \ \text{and}\  p\geq2, \text{and}\  n\geq4,
\end{cases}
\end{equation*}
then $M$ is isometric to a sphere form.
\end{thm}
\begin{rem}  When $n=3$ and $p\geq2$, or $n\geq 4$, the inequality  of this theorem is optimal.
The critical case is given by the following example. If $(\mathbb{S}^1(t)\times \mathbb{S}^{n-1}, g_t)$ is the product of the circle
of radius $t$ with $\mathbb{S}^{n-1}$, and if $g_t$ is the standard product metric normalized such that $Vol(g_t)=1$, we have
$W=0$, $g_t$ is a Yamabe metric for small $t$ (see \cite{S1}), and $\left(\int_M|\mathring{Ric}|^{p}\right)^{\frac 1p}=\frac{R(g_t)}{\sqrt{n(n-1)}}$, which is the critical case  of the inequality  in Theorem 1.13. We know that $(\mathbb{S}^1(t)\times \mathbb{S}^{n-1}, g_t)$ is not Einstein.

When $n\geq 5$ and $p=\frac n2$, Theorem 1.13 reduces to Thereom 1.1 in \cite{{HV}}. Theorem 1.13  improves Thereom 1.1 and Corollary 1 in \cite{{HV}}.
\end{rem}

As we mentioned above, Theorem 1.13 is sharped. By this we mean that we can precisely characterize the case of equality:
\begin{thm}
Let $M^n(n\geq 4)$ be a compact  locally conformally flat
Riemannian $n$-manifold with constant positive scalar curvature. For $p\geq \frac n2$,
if
\begin{eqnarray*}\left(\int_{M}|\mathring{Ric}|^{p}\right)^{\frac 1p}=C_3,\end{eqnarray*}
then i) $M$ is covered isometrically by $\mathbb{S}^1\times \mathbb{S}^{n-1}$ with the product metric;

ii) $M$ is covered isometrically by $(\mathbb{S}^1\times \mathbb{S}^{n-1}, dt^2+F^2(t)g_{\mathbb{S}^{n-1}})$, where $(\mathbb{S}^{n-1}, g_{\mathbb{S}^{n-1}})$ is a
round sphere and $F$ is a non-constant, positive,
periodic function satisfying a precise ODE. This metric is called a rotationally symmetric Derdzi\'{n}ski metric in \cite{{Ca},{D}}.
\end{thm}


\section{Proof of  Lemmas}
\label{}
In what follows, we adopt, without further comment, the moving frame notation with respect to a chosen local orthonormal frame.

Let $M$ be a  Riemannian manifold with harmonic curvature. The
decomposition of the Riemannian curvature tensor  into irreducible components yield
\begin{eqnarray*}
R_{ijkl}&=&W_{ijkl}+\frac{1}{n-2}(R_{ik}\delta_{jl}-R_{il}\delta_{jk}+R_{jl}\delta_{ik}-R_{jk}\delta_{il})\nonumber\\
&&-\frac{R}{(n-1)(n-2)}(\delta_{ik}\delta_{jl}-\delta_{il}\delta_{jk})\nonumber\\
&=&W_{ijkl}+\frac{1}{n-2}(\mathring{R}_{ik}\delta_{jl}-\mathring{R}_{il}\delta_{jk}+\mathring{R}_{jl}\delta_{ik}-\mathring{R}_{jk}\delta_{il})\nonumber\\
&&+\frac{R}{n(n-1)}(\delta_{ik}\delta_{jl}-\delta_{il}\delta_{jk}),
\end{eqnarray*}
where $R_{ijkl}$, $W_{ijkl}$, $R_{ij}$  and $\mathring{R}_{ij}$  denote the components of $Rm$, the Weyl curvature tensor $W$,  the Ricci tensor $Ric$ and the trace-free Ricci tensor $\mathring{Ric}=Ric-\frac{R}{n}g$, respectively,  and  $R$  is the scalar curvature.

 The trace-free Riemannian curvature tensor $\mathring{Rm}$ is
\begin{eqnarray}\mathring{R}_{ijkl}={R}_{ijkl}-\frac{R}{n(n-1)}(\delta_{ik}\delta_{jl}-\delta_{il}\delta_{jk}).\end{eqnarray}
Then the following equalities are easily obtained
from the properties of curvature tensor:
\begin{eqnarray}g^{ik}\mathring{R}_{ijkl}=\mathring{R}_{jl},\end{eqnarray}
\begin{eqnarray}\mathring{R}_{ijkl}+\mathring{R}_{iljk}+\mathring{R}_{iklj}=0,\end{eqnarray}
\begin{eqnarray}\mathring{R}_{ijkl}=\mathring{R}_{klij}=-\mathring{R}_{jikl}=-\mathring{R}_{ijlk},\end{eqnarray}
\begin{eqnarray}|\mathring{Rm}|^2=|W|^2+\frac{4}{n-2}|\mathring{Ric}|^2.\end{eqnarray}
Moreover, by the assumption of harmonic curvature,
we compute
\begin{eqnarray}
\mathring{R}_{ijkl,m}+\mathring{R}_{ijmk,l}+\mathring{R}_{ijlm,k}=0,
\end{eqnarray}
and
\begin{eqnarray}
\mathring{R}_{ijkl,l}=0.
\end{eqnarray}

Now, we compute the Laplacian of $|\mathring{Rm}|^2$.
\begin{lem}
Let $M$  be a complete  Riemannian $n$-manifold
with harmonic curvature.   Then
\begin{eqnarray}
\triangle|\mathring{Rm}|^2\geq2|\nabla \mathring{Rm}|^2-2C(n)|\mathring{Rm}|^3+2AR|\mathring{Rm}|^2,
\end{eqnarray}
where \begin{equation*}A=
\begin{cases}\frac{1}{n-1}, R\geq0\\ \frac{2}{n}, R<0,
\end{cases}\end{equation*}and $C(n)=2[\frac{2(n^2+n-4)}{\sqrt{(n-1)n(n+1)(n+2)}}+ \frac{n^2-n-4}{2\sqrt{(n-2)(n-1)n(n+1)}}+\sqrt{\frac{(n-2)(n-1)}{4n}}]$.
\end{lem}
\begin{rem} Although Lemma 2.1 has been proved in \cite{Ch}, we give an explicit coefficient of the term $|\mathring{Rm}|^3$ in (8).
When $M$ is a  complete locally
conformally flat Riemannian $n$-manifold, it follows from (10) that
\begin{equation*} \triangle|\mathring{Ric}|^2\geq2|\nabla \mathring{Ric}|^2-2\frac{n}{\sqrt{n(n-1)}}|\mathring{Ric}|^3+2\frac{R}{n-1}|\mathring{Ric}|^2.\end{equation*}
By the Kato inequality $|\nabla \mathring{Ric}|^2\geq \frac{n+2}{n}|\nabla |\mathring{Ric}||^2$, we obtain ( see \cite{{PRS},{XZ}})
\begin{equation} |\mathring{Ric}|\triangle|\mathring{Ric}|\geq\frac{2}{n}|\nabla |\mathring{Ric}||^2-\frac{n}{\sqrt{n(n-1)}}|\mathring{Ric}|^3+\frac{R}{n-1}|\mathring{Ric}|^2.\end{equation}
\end{rem}
\begin{proof}
By the Ricci identities, we obtain from (1)-(7)
\begin{eqnarray}
\triangle|\mathring{Rm}|^2&=&2|\nabla \mathring{Rm}|^2+2\langle \mathring{Rm}, \triangle \mathring{Rm}\rangle=2|\nabla \mathring{Rm}|^2+2\mathring{R}_{ijkl}\mathring{R}_{ijkl,mm}\nonumber\\
&=&2|\nabla \mathring{Rm}|^2+2\mathring{R}_{ijkl}(\mathring{R}_{ijkm,lm}+\mathring{R}_{ijml,km})\nonumber\\
&=&2|\nabla \mathring{Rm}|^2+4\mathring{R}_{ijkl}\mathring{R}_{ijkm,lm}\nonumber\\
&=&2|\nabla \mathring{Rm}|^2+4\mathring{R}_{ijkl}(\mathring{R}_{ijkm,ml}
+\mathring{R}_{hjkm}R_{hilm}\nonumber\\
&&+\mathring{R}_{ihkm}R_{hjlm}+\mathring{R}_{ijhm}R_{hklm}
+\mathring{R}_{ijkh}R_{hmlm})\nonumber\\
&=&2|\nabla \mathring{Rm}|^2+4\mathring{R}_{ijkl}(\mathring{R}_{hjkm}R_{hilm}
+\mathring{R}_{ihkm}R_{hjlm}\nonumber\\
&&+\mathring{R}_{ijhm}R_{hklm}
+\mathring{R}_{ijkh}R_{hmlm})\nonumber\\
&=&2|\nabla \mathring{Rm}|^2+4\mathring{R}_{ijkl}(\mathring{R}_{hjkm}\mathring{R}_{hilm}
+\mathring{R}_{ihkm}\mathring{R}_{hjlm}
+\mathring{R}_{ijhm}\mathring{R}_{hklm}
\nonumber\\&&+\mathring{R}_{ijkh}\mathring{R}_{hmlm})
+\frac{4R}{n(n-1)}\mathring{R}_{ijkl}(\mathring{R}_{ljki}+\mathring{R}_{ilkj}+\mathring{R}_{ijlk}\nonumber\\
&&+\mathring{R}_{jk}\delta_{il}-\mathring{R}_{ik}\delta_{jl})+\frac{4R}{n}|\mathring{Rm}|^2\nonumber\\
&=&2|\nabla \mathring{Rm}|^2+4\mathring{R}_{ijkl}(\mathring{R}_{hjkm}\mathring{R}_{hilm}
+\mathring{R}_{ihkm}\mathring{R}_{hjlm}
+\mathring{R}_{ijhm}\mathring{R}_{hklm}\nonumber\\
&&+\mathring{R}_{ijkh}\mathring{R}_{hl})
-\frac{8R}{n(n-1)}|\mathring{Ric}|^2+\frac{4R}{n}|\mathring{Rm}|^2\nonumber\\
&=&2|\nabla \mathring{Rm}|^2-4\mathring{R}_{ijlk}(2\mathring{R}_{jhkm}\mathring{R}_{himl}
+\frac 12\mathring{R}_{hmij}\mathring{R}_{lkhm}\nonumber\\
&&+\mathring{R}_{ijkh}\mathring{R}_{hl})
-\frac{8R}{n(n-1)}|\mathring{Ric}|^2+\frac{4R}{n}|\mathring{Rm}|^2.
\end{eqnarray}

We consider $\mathring{Rm}$ as a self adjoint operator on $\wedge^2 V$ and $S^2 V$. By the algebraic inequality for $m$-trace-free symmetric two-tensors $T$, i.e., $tr(T^3)\leq\frac{m-2}{\sqrt{m(m-1)}}|T|^3$, and the eigenvalues $\lambda_i$ of $T$ satisfy $|\lambda_i|\leq\sqrt{\frac{m-1}{m}}|T|$ in \cite{Hu}, we obtain
\begin{eqnarray}|\mathring{R}_{ijlk}(2\mathring{R}_{jhkm}\mathring{R}_{himl}
+\frac 12\mathring{R}_{hmij}\mathring{R}_{lkhm})|\leq2|\mathring{R}_{ijlk}\mathring{R}_{jhkm}\mathring{R}_{himl}|
+\frac 12|\mathring{R}_{ijlk}\mathring{R}_{hmij}\mathring{R}_{lkhm}|\nonumber\\
\leq[\frac{2(n^2+n-4)}{\sqrt{(n-1)n(n+1)(n+2)}}+ \frac{n^2-n-4}{2\sqrt{(n-2)(n-1)n(n+1)}}]|\mathring{Rm}|^3,
\end{eqnarray}
and
\begin{eqnarray}|\mathring{R}_{ijkl}\mathring{R}_{ijkh}\mathring{R}_{hl}|\leq\sqrt{\frac{n-1}{n}}|\mathring{Ric}||\mathring{Rm}|^2.
\end{eqnarray}
From (5), we have
\begin{equation}|\mathring{Ric}|^2\leq \frac{n-2}{4}|\mathring{Rm}|^2.
\end{equation}
Combining  with (10)-(13), we obtain  that
\begin{eqnarray*}
\triangle|\mathring{Rm}|^2\geq 2|\nabla \mathring{Rm}|^2+2AR|\mathring{Rm}|^2-4[\sqrt{\frac{(n-2)(n-1)}{4n}}
\\+\frac{2(n^2+n-4)}{\sqrt{(n-1)n(n+1)(n+2)}}+ \frac{n^2-n-4}{2\sqrt{(n-2)(n-1)n(n+1)}}]|\mathring{Rm}|^3.
\end{eqnarray*}
This completes the proof of this Lemma.
\end{proof}

\begin{lem}
Let $M$  be a complete noncompact Riemannian $n$-manifold
satisfying a Sobolev inequality of the following form:
\begin{eqnarray}
\left(\int_M|f|^{\frac{2n}{n-2}}\right)^{\frac{n-2}{n}}\leq D_1\int_M|\nabla f|^2+F_1\int_M|f|^2,  \forall f\in C^{\infty}_0(M).
\end{eqnarray}
 If a non-negative function $u \in C^{\infty}(M)$ satisfies $\int_M u^{\frac n2}<+\infty$ and
 \begin{eqnarray}
\triangle u\geq au^2+bu
\end{eqnarray}
for some constants $a$ and $b$.
 Then,
 given any $\epsilon>0$ and any $x_0\in M$ there exists a
geodesic ball $B_{r}(x_0)$ with center $x_0$ and radius $r$ such that $u(x)<\epsilon$ for all
$x\in M\setminus B_{r}(x_0)$.
\end{lem}
\begin{rem} We can carry out the proof of Lemma 2.3  by suitable modification to  the proof of Theorem 1.1 in \cite{BCW}.
\end{rem}
\begin{proof} Let us fix a point $x_0\in M$, consider the following open domains:  $$E(T)=\{x\in M|d(x_0, x)>T\}, A(T, S)=\{x\in M|T<d(x_0, x)<S\}.$$

First, take a cut-off function $\varphi \in C^{\infty}_0(M)$ with the properties $Supp(\varphi)\subset A(r, s+2)\subset E(r), r+2<s$ and $\varphi=1$ on $A(r+2, s)$. Using Cauchy-Schwarz inequality,  multiplying  (15) by $\varphi^2u^{\frac n2-1}$ and integrating,
 we have
\begin{eqnarray*}
\frac{8(n-2)}{n^2}\int_M\varphi^2|\nabla u^{\frac{n}{4}}|^2\leq-a\int_M\varphi^2u^{\frac n2+1}-b\int_M\varphi^2u^{\frac n2}\\+\frac{4}{n^2}\int_M\varphi^2|\nabla u^{\frac{n}{4}}|+4\int_Mu^{\frac{n}{2}}|\nabla\varphi|^2,
\end{eqnarray*}
which gives
\begin{eqnarray}
\int_M\varphi^2|\nabla u^{\frac{n}{4}}|^2\leq-na\int_M\varphi^2u^{\frac n2+1}-n(b-4)\int_{Supp(\varphi)}\varphi^2u^{\frac n2}.
\end{eqnarray}
We now apply Sobolev inequality (14) to the function $f=\varphi u^{\frac{n}{4}}$ and  obtain
\begin{eqnarray}
\left(\int_M|\varphi u^{\frac{n}{4}}|^{\frac{2n}{n-2}}\right)^{\frac{n-2}{n}}\leq 2D_1\int_M\varphi^2|\nabla u^{\frac{n}{4}}|^2+2D_1\int_Mu^{\frac{n}{2}}|\nabla \varphi|^2+F_1\int_M\varphi^2 u^{\frac{n}{2}}.
\end{eqnarray}
Plugging (16) into (17), we obtain
\begin{eqnarray}
\left(\int_M\varphi^{\frac{2n}{n-2}} u^{\frac{n^2}{2(n-2)}}\right)^{\frac{n-2}{n}}\leq -2naD_1\int_M\varphi^2u^{\frac n2+1}+F_2\int_{Supp(\varphi)}u^{\frac n2}.
\end{eqnarray}
Using H\"{o}lder inequality, we find that
\begin{equation}
\int_M\varphi^2u^{\frac n2+1}\leq\left(\int_{Supp(\varphi)}u^{\frac n2}\right)^{\frac{2}{n}}\left(\int_M\varphi^{\frac{2n}{n-2}} u^{\frac{n^2}{2(n-2)}}\right)^{\frac{n-2}{n}}
\end{equation}
and we deduce from (18) that
\begin{eqnarray}
\left(1+2naD_1\left(\int_{Supp(\varphi)}u^{\frac n2}\right)^{\frac{2}{n}}\right)\left(\int_M\varphi^{\frac{2n}{n-2}} u^{\frac{n^2}{2(n-2)}}\right)^{\frac{n-2}{n}}\leq F_2\int_{Supp(\varphi)}u^{\frac n2}.
\end{eqnarray}
If we choose $r$ sufficiently big such that $\left(1+2naD_1\left(\int_{Supp(\varphi)}u^{\frac n2}\right)^{\frac{2}{n}}\right)\geq\frac12$, then we get
\begin{eqnarray}
\left(\int_{E(r+2)} u^{\frac{n^2}{2(n-2)}}\right)^{\frac{2(n-2)}{n^2}}\leq F_3\left(\int_{E(r)}u^{\frac n2}\right)^{\frac{2}{n}}.
\end{eqnarray}

Second, we multiply (15) by $\varphi^2u^{k-1}$, for $k\geq \frac 32$,  and integrate to find
\begin{eqnarray*}
\frac{4(k-1)}{k^2}\int_M\varphi^2|\nabla u^{\frac{k}{2}}|^2\leq-a\int_M\varphi^2u^{k+1}-b\int_M\varphi^2u^{k}+\frac{4}{k}\int_Mu^{\frac{k}{2}}|\nabla\varphi|\varphi|\nabla u^{\frac{k}{2}}|,
\end{eqnarray*}
which gives
\begin{eqnarray}
\frac{4k-5}{k^2}\int_M\varphi^2|\nabla u^{\frac{k}{2}}|^2\leq-a\int_M\varphi^2u^{k+1}-b\int_M\varphi^2u^{k}+4\int_Mu^{k}|\nabla\varphi|^2.
\end{eqnarray}
Applying Sobolev inequality (14) to the function $f=\varphi u^{\frac{k}{2}}$, we have
\begin{eqnarray}
\left(\int_M|\varphi u^{\frac{k}{2}}|^{\frac{2n}{n-2}}\right)^{\frac{n-2}{n}}\leq 2D_1\int_M\varphi^2|\nabla u^{\frac{k}{2}}|^2+2D_1\int_Mu^{k}|\nabla \varphi|^2+F_1\int_M|\varphi u^{\frac{k}{2}}|^2.
\end{eqnarray}
Plugging (22) into (23), we get
\begin{eqnarray}
\left(\int_M|\varphi u^{\frac{k}{2}}|^{\frac{2n}{n-2}}\right)^{\frac{n-2}{n}}\leq D_2k(\int_M\varphi^2u^{k+1}+\int_Mu^{k}|\nabla \varphi|^2+\int_M\varphi^2 u^{k}).
\end{eqnarray}
Using H\"{o}lder inequality, from (24) we get
\begin{eqnarray*}
\left(\int_M|\varphi u^{\frac{k}{2}}|^{\frac{2n}{n-2}}\right)^{\frac{n-2}{n}}\leq D_2k[\left(\int_M|\varphi u^{\frac{k}{2}}|^{\frac{2t}{t-2}}\right)^{\frac{t-2}{t}}\left(\int_{Supp(\varphi)}u^{\frac{t}{2}}\right)^{\frac{2}{t}}\\+\int_Mu^{k}|\nabla \varphi|^2+\int_M\varphi^2 u^{k}],
\end{eqnarray*}
where $t=\frac{n^2}{n-2}$.
By (21) and the above inequality, we conclude that
\begin{eqnarray}
\left(\int_M|\varphi u^{\frac{k}{2}}|^{\frac{2n}{n-2}}\right)^{\frac{n-2}{n}}\leq D_3k\left[\left(\int_M|\varphi u^{\frac{k}{2}}|^{\frac{2t}{t-2}}\right)^{\frac{t-2}{t}}+\int_Mu^{k}|\nabla \varphi|^2+\int_M\varphi^2 u^{k}\right].
\end{eqnarray}
Since $2<\frac{2t}{t-2}<\frac{2n}{n-2}$, we have,
by interpolation,
\begin{eqnarray*}
\left(\int_M|\varphi u^{\frac{k}{2}}|^{\frac{2t}{t-2}}\right)^{\frac{t-2}{t}}\leq2\epsilon\left(\int_M|\varphi u^{\frac{k}{2}}|^{\frac{2n}{n-2}}\right)^{\frac{n-2}{n}}+2\epsilon^{-\frac{n-2}{2}}\int_M\varphi^2 u^{k}.
\end{eqnarray*}
We choose $\epsilon$ such that $4D_3k\epsilon=1$ and  plug the above inequality into (25) to
obtain
\begin{eqnarray}
\left(\int_M|\varphi u^{\frac{k}{2}}|^{\frac{2n}{n-2}}\right)^{\frac{n-2}{n}}\leq D_4k\left[\left(1+k^{\frac{n-2}{n}}\right)\int_M\varphi^2 u^{k}+\int_Mu^{k}|\nabla \varphi|^2\right].
\end{eqnarray}

Finally, we now perform the iteration method under (26).
 Define $T_0=T, T_{i+1}=T_i+2^{-(i+1)}T$, and $S_0=2S, S_{i+1}=S_i-2^{-(i+1)}S$.
Choose $\varphi_i \in C^{\infty}_0(M)$ such that $0\leq\varphi_i\leq1, |\nabla \varphi_i|\leq5^i,\sup\varphi_i\subset A_i\triangleq A(T_i, S_i)$ and $\varphi_i\mid A_{i+1}=1$. let $\chi_i=\textbf{1}_{\sup\varphi_i}$. By (26), we get
\begin{eqnarray*}
\left(\int_M|\varphi_i u^{\frac{k}{2}}|^{\frac{2n}{n-2}}\right)^{\frac{n-2}{n}}&\leq& D_4k\left[\left(1+k^{\alpha\beta}\right)\int_M\chi_i^2 u^{k}+25^i\int_Mu^{k}\chi_i^2\right]\nonumber\\
&\leq&D_4k(26^i+k^{\alpha\beta})\int_Mu^{k}\chi_i^2,
\end{eqnarray*}
where $\beta=\frac{(n-2)^2}{2n}$. We can rewrite the above inequality as
\begin{eqnarray*}
\left(\int_{A_{i+1}}u^{k\alpha}\right)^{\frac{1}{\alpha}}\leq D_4k(26^i+k^{\alpha\beta})\int_{A_{i}}u^{k}.
\end{eqnarray*}
Take $k=k_i=\gamma\alpha^i$, for $i\geq0$, then with obvious notations
\begin{eqnarray*}
\left(\int_{A_{i+1}}u^{k_{i+1}}\right)^{\frac{1}{k_{i+1}}}&\leq& \{D_4\gamma\alpha^i[26^i+(\gamma\alpha^i)^{\alpha\beta}]\}^{\frac{1}{\gamma\alpha^i}}(\int_{A_{i}}u^{k_i})^{\frac{1}{k_{i}}}\nonumber\\
&\leq&[D_4\gamma(1+\gamma^{\alpha\beta})]^{\frac{1}{\gamma\alpha^i}}\{(\alpha^{\frac{1}{\gamma}}\max\{26,\alpha^{\alpha\beta}\})^i\}^{\frac{1}{\alpha^i}}
(\int_{A_{i}}u^{k_i})^{\frac{1}{k_{i}}}\nonumber\\
&\leq& D_5e^{iD_6\alpha^{-i}}(\int_{A_{i}}u^{k_i})^{\frac{1}{k_{i}}}.
\end{eqnarray*}
Therefore,
\begin{eqnarray*}
\|u\|_{\infty, A_\infty}=\|u\|_{\infty, A(2T, S)}\leq D_5e^{D_6\sum_{i=0}^{\infty}i\alpha^{-i}}(\int_{A(T, 2S)}u^{\gamma})^{\frac{1}{\gamma}}
\end{eqnarray*}
Taking $\gamma=\frac n2$, we conclude that $\|u\|_{\infty, A(2T, S)}\leq A_M \|u\|_{n, E(T)}$ and the
assertion in Lemma 2.3 follows by letting $S$ tend to infinity.
\end{proof}

\section{Proof of Theorems}
\label{}
\begin{proof}[{\bf Proof of Theorem
1.1}]
From (8), by the Kato inequality $|\nabla \mathring{Rm}|^2\geq |\nabla |\mathring{Rm}||^2$, we obtain
\begin{eqnarray}
|\mathring{Rm}|\triangle|\mathring{Rm}|=\frac12\triangle|\mathring{Rm}|^2-|\nabla |\mathring{Rm}||^2\geq-C(n)|\mathring{Rm}|^3+AR|\mathring{Rm}|^2.
\end{eqnarray}
Let $u=|\mathring{Rm}|$.
By (27), we compute
\begin{eqnarray} u^{\alpha}\triangle u^{\alpha}&=&u^{\alpha}\left(\alpha(\alpha-1)u^{\alpha-2}|\nabla u|^2+\alpha u^{\alpha-1}\triangle u\right)\nonumber\\
&=&\frac{\alpha-1}{\alpha}|\nabla u^{\alpha}|^2+\alpha
u^{2\alpha-2}u\triangle u\nonumber\\
&\geq&\frac{\alpha-1}{\alpha}|\nabla u^{\alpha}|^2-C(n)\alpha u^{2\alpha+1}+\alpha AR u^{2\alpha},
\end{eqnarray}
where $\alpha$ is a positive constant. Taking $\alpha=\frac{2p}{n}\geq 1$, using the Young's inequality, from (28) we obtain
\begin{eqnarray}
u^{\alpha}\triangle u^{\alpha}\geq au^{3\alpha}+b u^{2\alpha}.
\end{eqnarray}
where $a$ and $b$ are two constants depending only on $n, \alpha$ and $R$. Setting $w=u^{\alpha}$, we can rewrite (29) as
\begin{eqnarray}
\triangle w\geq aw^2+b w.
\end{eqnarray}
Since $M$ has the positive Yamabe constant or satisfies the Sobolev inequality, combining with (30), by Lemma 2.3, we can prove Theorem 1.1.
\end{proof}
\begin{proof}[{\bf Proof of Theorem
1.2}]
By (1), we have
\begin{eqnarray}{R}_{ijij}=\mathring{R}_{ijij}+\frac{R}{n(n-1)}.\end{eqnarray}
Note that $R$
is positive. From (31), we see from Theorem 1.1 that there is a positive constant $\delta$ such that ${R}_{ijij}>\delta$ in $M\setminus\Omega$
for some compact set $\Omega$. This implies
that the Ricci curvature is bounded from below by a positive constant outside some
geodesic sphere, hence the manifold is compact (for detail, see Lemma 3.5 of \cite{XZ}).
\end{proof}
\begin{proof}[{\bf Proof of Theorem
1.4}] When $R>0$, we see from Theorem 1.2 that $M$ is compact. Taking $\alpha=\frac{2p}{n}\geq1$.
From (28), using the Young's inequality, we have
\begin{eqnarray}
u^{\alpha}\triangle u^{\alpha}\geq\frac{\alpha-1}{\alpha}|\nabla u^{\alpha}|^2
-C(n)\epsilon^{1-\alpha}u^{3\alpha}-[C(n) (\alpha-1)\epsilon-\alpha A R]u^{2\alpha}.
\end{eqnarray}
Setting $w=u^{\alpha}$, we can rewrite (32) as
\begin{eqnarray}
w\triangle w\geq\frac{\alpha-1}{\alpha}|\nabla w|^2
-C(n)\epsilon^{1-\alpha}w^{3}-[C(n) (\alpha-1)\epsilon-\alpha A R]w^{2}.
\end{eqnarray}
From (33), we derive
\begin{eqnarray}
w^{\beta}\triangle w^{\beta}
\geq(1-\frac{1}{\alpha\beta})|\nabla w^{\beta}|^2
-C(n)\beta\epsilon^{1-\alpha}w^{2\beta+1}\nonumber\\-\beta[C(n) (\alpha-1)\epsilon-\alpha A R]w^{2\beta},
\end{eqnarray}
where $\beta$ is a positive constant. Integrating by parts (34), we get
\begin{eqnarray}
(2-\frac{1}{\alpha\beta})\int_{M}|\nabla w^\beta|^2-C(n)\beta\epsilon^{1-\alpha}\int_{M}w^{2\beta+1}\nonumber\\
-\beta[C(n) (\alpha-1)\epsilon-\frac{R\alpha}{n-1}]\int_{M}w^{2\beta}\leq 0.
\end{eqnarray}
From (35), using the H\"{o}lder inequality, we have
\begin{eqnarray}
(2-\frac{1}{\alpha\beta})\int_{M}|\nabla w^\beta|^2-C(n)\beta\epsilon^{1-\alpha}(\int_{M}w^{\frac{2n\beta}{n-2}})^{\frac{n-2}{n}}(\int_{M}w^{\frac n2})^{\frac{2}{n}}\nonumber\\
-\beta[C(n) (\alpha-1)\epsilon-\frac{R\alpha}{n-1}]\int_{M}w^{2\beta}\leq 0.
\end{eqnarray}

Case 1. When $3\leq n\leq5$ and $1\leq\alpha<\frac{4}{n-2}$, if $\alpha>1$, set $\epsilon=\frac{(6-n)\alpha R}{4(n-1)(\alpha-1)C(n)}$; if $\alpha=1$, set $\epsilon=1$. Take $\beta=\frac 1\alpha$. By the definition of Yamabe constant $Q(M)$, from (36) we get
\begin{eqnarray}
\left[Q(M)-\frac{C(n)\epsilon^{1-\alpha}}{\alpha}\left(\int_{M}|\mathring{Rm}|^{p}\right)^{\frac{2}{n}}\right]
\left(\int_{M}w^{\frac{2n\beta}{n-2}}\right)^{\frac{n-2}{n}}\leq 0.
\end{eqnarray}
We choose $\left(\int_{M}|\mathring{Rm}|^{p}\right)^{\frac 1p}<C_1$ such that (37) implies $\left(\int_{M}w^{\frac{2n\beta}{n-2}}\right)^{\frac{n-2}{n}}=0$, that is, $|\mathring{Rm}|=0$, i.e., $M$ is  Einstein manifold and locally conformally flat manifold. Hence $M$ is isometric to a spherical
space form.

Case 2. When $3\leq n\leq5$ and $\alpha\geq\frac{4}{n-2}$, and $n\geq6$, set $\epsilon=\frac{R}{(n-1)C(n)}$ and
 $\frac{1}{\alpha\beta}=1+\sqrt{1-\frac{4}{(n-2)\alpha}}$. We also get
\begin{eqnarray}
\left[(2-\frac{1}{\alpha\beta})Q(M)-C(n)\beta\epsilon^{1-\alpha}\left(\int_{M}|\mathring{Rm}|^{p}\right)^{\frac{2}{n}}\right]
\left(\int_{M}w^{\frac{2n\beta}{n-2}}\right)^{\frac{n-2}{n}}\leq 0.
\end{eqnarray}
We choose $\left(\int_{M}|\mathring{Rm}|^{p}\right)^{\frac 1p}<C_1$ such that (38) implies $\left(\int_{M}w^{\frac{2n\beta}{n-2}}\right)^{\frac{n-2}{n}}=0$, that is, $|\mathring{Rm}|=0$, i.e., $M$ is  Einstein manifold and locally conformally flat manifold. Hence $M$ is isometric to a spherical
space form.
\end{proof}

\begin{proof}[{\bf Proof of Theorem
1.6}]
Let $\phi$ be a smooth compactly supported function on $M$.  Multiplying (34) by
$\phi^2$ and integrating over $M$, we obtain
\begin{eqnarray*} (1-\frac{1}{\alpha\beta})\int_{M}|\nabla w^{\beta}|^2 \phi^2\leq {C(n)\beta\epsilon^{1-\alpha}}\int_{M}w^{2\beta+1}\phi^2+\int_{M}w^{\beta}\phi^2\triangle
w^{\beta}\\
+\beta[C(n) (\alpha-1)\epsilon-\frac{2R\alpha}{n}]\int_{M}w^{2\beta}\phi^2\\
=C(n)\beta\epsilon^{1-\alpha}\int_{M}w^{2\beta+1}\phi^2-2\int_{M}w^{\beta}\phi\langle\nabla \phi,\nabla
w^{\beta}\rangle\\
-\int_{M}|\nabla
w^{\beta}|^2
\phi^2+\beta[C(n) (\alpha-1)\epsilon-\frac{2R\alpha}{n}]\int_{M}w^{2\beta}\phi^2,
\end{eqnarray*}
which gives
\begin{eqnarray}(2-\frac{1}{\alpha\beta})\int_{M}|\nabla w^{\beta}|^2 \phi^2\leq C(n)\beta\epsilon^{1-\alpha}\int_{M}w^{2\beta+1}\phi^2-2\int_{M}w^{\beta}\phi\langle\nabla \phi,\nabla
w^{\beta}\rangle\nonumber\\
+\beta[C(n) (\alpha-1)\epsilon-\frac{2R\alpha}{n}]\int_{M}w^{2\beta}\phi^2.
\end{eqnarray}
Using the Cauchy-Schwarz inequality, we can rewrite (39) as
\begin{eqnarray} \left(2-\frac{1}{\alpha\beta}-\varepsilon\right)\int_{M}|\nabla w^{\beta}|^2 \phi^2\leq C(n)\beta\epsilon^{1-\alpha}\int_{M}w^{2\beta+1}\phi^2+\frac{1}{\varepsilon}\int_M
w^{2\beta}|\nabla \phi|^2\nonumber\\
+\beta[C(n) (\alpha-1)\epsilon-\frac{2R\alpha}{n}]\int_{M}w^{2\beta}\phi^2,
\end{eqnarray}
for the positive constant $\varepsilon$. By the definition of Yamabe constant $Q(M)$ and (40),  we have
\begin{eqnarray} Q(M)\left(\int_{M}(\phi w^{\beta})^{\frac{2n}{n-2}}\right)^{\frac{n-2}{n}}\leq \int_{M}\left(|\nabla (\phi w^{\beta})|^2+\frac{(n-2)R
w^{2\beta} \phi^2}{4(n-1)}\right)\nonumber\\
=\int_{M}(w^{2\beta}|\nabla \phi |^2+\phi^2|\nabla w^{\beta}|^2+2\phi w^{\beta}\langle\nabla \phi,\nabla
w^{\beta}\rangle+\frac{(n-2)R
w^{2\beta} \phi^2}{4(n-1)})\nonumber\\
\leq(1+\frac{1}{\eta})\int_{M}w^{2\beta}|\nabla \phi |^2+(1+\eta)\int_{M}\phi^2|\nabla w^{\beta}|^2
+\int_{M}\frac{(n-2)R
w^{2\beta} \phi^2}{4(n-1)}\nonumber\\
\leq G\int_{M}w^{2\beta}|\nabla \phi |^2+H\int_{M}w^{2\beta+1}\phi^2
+I\int_{M}
w^{2\beta} \phi^2,
\end{eqnarray}
where
\begin{eqnarray*}G&=&1+\frac{1}{\eta}+\frac{1+\eta}{\varepsilon(2-\frac{1}{\alpha\beta}-\varepsilon)},\\
H&=&\frac{(1+\eta)C(n)\beta\epsilon^{1-\alpha}}{(2-\frac{1}{\alpha\beta}-\varepsilon)},\\
I&=&\frac{(n-2)R}{4(n-1)}-\frac{2(1+\eta)\alpha\beta R}{n(2-\frac{1}{\alpha\beta}-\varepsilon)}+\frac{(1+\eta)C(n)(\alpha-1)\beta\epsilon}{(2-\frac{1}{\alpha\beta}-\varepsilon)}.\end{eqnarray*}

We first consider the case of $\gamma\in (1, \frac{n(n-2)+\sqrt{n(n-2)(n^2-10n+8)}}{4(n-1)})$.
When $n\geq 10$, noting that $\epsilon$, $\varepsilon$ and $\eta$ are sufficiently small, we choose $\frac12<\alpha\beta<\frac{n(n-2)+\sqrt{n(n-2)(n^2-10n+8)}}{8(n-1)}$
such that $I\leq0$. Thus from (41) we have
\begin{eqnarray*} Q(M)\left(\int_{M}(\phi w^{\beta})^{\frac{2n}{n-2}}\right)^{\frac{n-2}{n}}
\leq G\int_{M}w^{2\beta}|\nabla \phi |^2+H\int_{M}w^{2\beta+1}\phi^2\nonumber\\
\leq G\int_{M}w^{2\beta}|\nabla \phi |^2+H\left(\int_{M}(\phi w^{\beta})^{\frac{2n}{n-2}}\right)^{\frac{n-2}{n}}\left(\int_{M} w^{\frac n2}\right)^{\frac{2}{n}}.
\end{eqnarray*}
Since $\int_{M} w^{\frac n2}=\int_{M} u^{p}$ is sufficiently small, the second term in
the right-hand side of the above can be absorbed in the left-hand side. Therefore,
there exists a constant $J>0$, such that
\begin{eqnarray} J\left(\int_{M}(\phi u^{\alpha\beta})^{\frac{2n}{n-2}}\right)^{\frac{n-2}{n}}
\leq G\int_{M}u^{2\alpha\beta}|\nabla \phi |^2.
\end{eqnarray}
Let us choose a cutoff function $\phi$
satisfying the properties that
\begin{equation*}\phi(x)=
\begin{cases}1 \quad\text{on}\quad B(r)\\0 \quad\text{on}\quad M\setminus
B(2r),
\end{cases}
\end{equation*}
and $|\nabla \phi|\leq \frac{2}{r}$. In particular,
if $M$ is compact, and if $r>d$, where $d$ is the diameter of $M$,  then $\phi=1$ on $M$. From (42), we get
\begin{eqnarray} J\left(\int_{B_{r}} u^{\frac{2n}{n-2}\alpha\beta}\right)^{\frac{n-2}{n}}
\leq \frac{4}{r^2}B\int_{M}u^{2\alpha\beta}.
\end{eqnarray}
Let $r\rightarrow +\infty$, by assumption that $\int_{M}u^{2\alpha\beta}<\infty$, from (43), we have $u=0$, i.e., $M$ is  Einstein manifold and locally conformally flat manifold. Hence  $M$ is isometric to a hyperbolic space form.

In the case of $\gamma\in (0, 1]$. Since $\int_{M}|\mathring{Rm}|^{p}<C,$ by Theorem 1.1, $|\mathring{Rm}|$ is bounded. Hence $\int_{M}|\mathring{Rm}|^{\gamma+1}<\infty$ for $\int_{M}|\mathring{Rm}|^{\gamma}<\infty$. For $\gamma+1\in(1, \frac{n(n-2)+\sqrt{n(n-2)(n^2-10n+8)}}{4(n-1)})$, we apply the above result to prove Theorem 1.6.
\end{proof}

\begin{proof}[{\bf Proof of Corollary
1.11}]Taking $\alpha=\frac{2p}{n}\geq1$. From (9), proceeding as in the proof of (34), we have

\begin{eqnarray}
w^{\beta}\triangle w^{\beta}
\geq(1-\frac{n-2}{n\alpha\beta})|\nabla w^{\beta}|^2
-\frac{n\beta\epsilon^{1-\alpha}}{\sqrt{n(n-1)}}w^{2\beta+1}\nonumber\\
-\beta[\frac{n(\alpha-1)\epsilon}{\sqrt{n(n-1)}}-\frac{R\alpha }{n-1}]w^{2\beta},
\end{eqnarray}
where $\beta$ is a positive constant. Based on (44), proceeding as in the proof of (40), we obtain
\begin{eqnarray} \left(2-\frac{n-2}{n\alpha\beta}-\varepsilon\right)\int_{M}|\nabla w^{\beta}|^2 \phi^2\leq \frac{n\beta\epsilon^{1-\alpha}}{\sqrt{n(n-1)}}\int_{M}w^{2\beta+1}\phi^2\nonumber\\+\frac{1}{\varepsilon}\int_M
w^{2\beta}|\nabla \phi|^2
+\beta[\frac{n(\alpha-1)\epsilon}{\sqrt{n(n-1)}}-\frac{R\alpha }{n-1}]\int_{M}w^{2\beta}\phi^2,
\end{eqnarray}
for the positive constant $\varepsilon$.

On the other hand, when $M$ is a  complete, simply connected, locally
conformally flat Riemannian n-manifold, $M$ satisfies the Sobolev inequality (see Corollary 3.2 in  \cite{H}):
\begin{equation*}\left(\int_{M} |f|^{\frac{2n}{n-2}}\right)^{\frac{n-2}{n}}\leq\frac{4}{n(n-2)\omega_{n}^{\frac2n}}\int_{M}\left(|\nabla f|^2+\frac{(n-2)}{4(n-1)}R
f^{2} \right), f\in C_0^{\infty}(M).\end{equation*}
Combining the above inequality with (45), we obtain
\begin{eqnarray}  \frac{n(n-2)\omega_{n}^{\frac2n}}{4}\left(\int_{M}(\phi w^{\beta})^{\frac{n}{n-2}}\right)^{\frac{n-2}{n}}\leq \int_{M}\left(|\nabla (\phi w^{\beta})|^2+\frac{(n-2)R
w^{2\beta} \phi^2}{4(n-1)}\right)\nonumber\\
=\int_{M}(w^{2\beta}|\nabla \phi |^2+\phi^2|\nabla w^{\beta}|^2+2\phi w^{\beta}\langle\nabla \phi,\nabla
w^{\beta}\rangle+\frac{(n-2)R
w^{2\beta} \phi^2}{4(n-1)})\nonumber\\
\leq(1+\frac{1}{\eta})\int_{M}w^{2\beta}|\nabla \phi |^2+(1+\eta)\int_{M}\phi^2|\nabla w^{\beta}|^2
+\int_{M}\frac{(n-2)R
w^{2\beta} \phi^2}{4(n-1)}\nonumber\\
\leq K\int_{M}w^{2\beta}|\nabla \phi |^2+L\int_{M}w^{2\beta+1}\phi^2
+M\int_{M}
w^{2\beta} \phi^2,
\end{eqnarray}
where
\begin{eqnarray*}K&=&1+\frac{1}{\eta}+\frac{1+\eta}{\varepsilon(2-\frac{(n-2)}{n\alpha\beta}-\varepsilon)},\\
L&=&\frac{(1+\eta)n\beta\epsilon^{1-\alpha}}{\sqrt{n(n-1)}(2-\frac{(n-2)}{n\alpha\beta}-\varepsilon)},\\
M&=&\frac{(n-2)R}{4(n-1)}-\frac{(1+\eta)\alpha \beta R}{(n-1)(2-\frac{(n-2)}{n\alpha\beta}-\varepsilon)}+\frac{(1+\eta)n(\alpha-1)\beta\epsilon}{\sqrt{n(n-1)}(2-\frac{(n-2)}{n\alpha\beta}-\varepsilon)}.\end{eqnarray*}
When $n\geq 6$, noting that $\epsilon$, $\varepsilon$ and $\eta$ are sufficiently small, we choose $\frac12\leq\alpha\beta<\frac{(n-2)(1+\sqrt{1-\frac 4n})}{4}$
such that $M\leq0$. Thus from (46) we have
\begin{eqnarray} \left(\frac{n(n-2)\omega_{n}^{\frac2n}}{4}-L\left(\int_{M} |\mathring{Rm}|^{p}\right)^{\frac{2}{n}}\right)\left(\int_{M}(\phi w^{\beta})^{\frac{2n}{n-2}}\right)^{\frac{n-2}{n}}
\leq K\int_{M}|\mathring{Rm}|^{2\alpha\beta}|\nabla \phi |^2.
\end{eqnarray}
So we choose $\left(\int_{M}|\mathring{Rm}|^{p}\right)^{\frac{2}{n}}$ small enough such that $N=\frac{n(n-2)\omega_{n}^{\frac2n}}{4}-L\left(\int_{M} |\mathring{Rm}|^{p}\right)^{\frac{2}{n}}>0$.

When $p=\frac n2$,  i.e., $\alpha=1$. We choose $\eta$
such that $M=0$, i.e., $(2-\frac{n-2}{n\alpha\beta}-\varepsilon)=\frac{4(1+\eta)\beta}{(n-2)}$. Thus we have
$$\frac{\frac{n(n-2)\omega_{n}^{\frac2n}}{4}}{L}=\sqrt{n(n-1)}\omega_{n}^{\frac2n}.$$
So we choose $\left(\int_{M}|\mathring{Rm}|^{\frac{n}{2}}\right)^{\frac{2}{n}}<\sqrt{n(n-1)}\omega_{n}^{\frac2n}$ such that $N=\frac{n(n-2)\omega_{n}^{\frac2n}}{4}-L\left(\int_{M} |\mathring{Rm}|^{\frac{n}{2}}\right)^{\frac{2}{n}}>0$.

Taking $\beta=\frac{n}{4}$. We rewrite (47) as
\begin{eqnarray*} N\left(\int_{M}(\phi w^{\beta})^{\frac{2n}{n-2}}\right)^{\frac{n-2}{n}}
\leq K\int_{M}|\mathring{Rm}|^{p}|\nabla \phi |^2.
\end{eqnarray*}
The rest of the proof runs as before.
Hence this completes the proof of Corollary
1.10.
\end{proof}
\begin{proof}[{\bf Proof of Theorem
1.13}]
Let $u=|\mathring{Ric}|$.
By (9), we compute
\begin{eqnarray} u^{\alpha}\triangle u^{\alpha}\geq(1-\frac{n-2}{n\alpha})|\nabla u^{\alpha}|^2
-\frac{n\alpha}{\sqrt{n(n-1)}}u^{2\alpha+1}+\frac{R\alpha }{n-1} u^{2\alpha},
\end{eqnarray}
where $\alpha$ is a positive constant.  Taking $\alpha=\frac{2p}{n}\geq1$.
From (48), using the Young's inequality, we have
\begin{eqnarray}
u^{\alpha}\triangle u^{\alpha}\geq(1-\frac{n-2}{n\alpha})|\nabla u^{\alpha}|^2
-\frac{n\epsilon^{1-\alpha}}{\sqrt{n(n-1)}}u^{3\alpha}-[\frac{n(\alpha-1)\epsilon}{\sqrt{n(n-1)}}-\frac{R\alpha }{n-1}]u^{2\alpha}.
\end{eqnarray}
Setting $w=u^{\alpha}$.
Based on (49), using the same argument as in the proof of (36), we obtain
\begin{eqnarray}
(2-\frac{n-2}{n\alpha\beta})\int_{M}|\nabla w^\beta|^2-\frac{n\beta\epsilon^{1-\alpha}}{\sqrt{n(n-1)}}(\int_{M}w^{\frac{2n\beta}{n-2}})^{\frac{n-2}{n}}(\int_{M}w^{\frac n2})^{\frac{2}{n}}\nonumber\\
-\beta[\frac{n(\alpha-1)\epsilon}{\sqrt{n(n-1)}}-\frac{R\alpha }{n-1}]\int_{M}w^{2\beta}\leq 0.
\end{eqnarray}

Case 1. When $n=3$ and $1\leq\alpha<\frac{4}{3}$, if $\alpha>1$, set $\epsilon=\frac{\sqrt{6}\alpha R}{24(\alpha-1)}$; if $\alpha=1$, set $\epsilon=1$. Take $\alpha\beta=\frac 13$. By the definition of Yamabe constant $Q(M)$, from (50) we get
\begin{eqnarray}
\left[Q(M)-\frac{3\epsilon^{1-\alpha}}{\sqrt{6}\alpha}\left(\int_{M}|\mathring{Ric}|^{p}\right)^{\frac{2}{n}}\right]
\left(\int_{M}w^{\frac{2n\beta}{n-2}}\right)^{\frac{n-2}{n}}\leq 0.
\end{eqnarray}
We choose $\left(\int_{M}|\mathring{Ric}|^{p}\right)^{\frac 1p}<C_3$ such that (50) implies $\left(\int_{M}w^{\frac{2n\beta}{n-2}}\right)^{\frac{n-2}{n}}=0$, that is, $|\mathring{Ric}|=0$, i.e., $M$ is  Einstein manifold.  Since $M$ is locally conformally flat manifold, $M$ is isometric to a spherical
space form.

Case 2. When $n=3$ and $\alpha\geq\frac{4}{3}$, and $n\geq4$, set $\epsilon=\frac{R}{\sqrt{(n-1)n}}$ and
 $\frac{1}{\alpha\beta}=\frac{n}{n-2}(1+\sqrt{1-\frac{4}{n\alpha}})$. We also get
\begin{eqnarray}
\left[(2-\frac{n-2}{n\alpha\beta})Q(M)-\frac{n\beta\epsilon^{1-\alpha}}{\sqrt{n(n-1)}}\left(\int_{M}|\mathring{Ric}|^{p}\right)^{\frac{2}{n}}\right]
\left(\int_{M}w^{\frac{2n\beta}{n-2}}\right)^{\frac{n-2}{n}}\leq 0.
\end{eqnarray}
We choose $\left(\int_{M}|\mathring{Ric}|^{p}\right)^{\frac 1p}<C_3$ such that (52) implies $\left(\int_{M}w^{\frac{2n\beta}{n-2}}\right)^{\frac{n-2}{n}}=0$, that is, $|\mathring{Ric}|=0$, i.e., $M$ is  Einstein manifold. Since $M$ is locally conformally flat manifold, $M$ is isometric to a spherical
space form.
\end{proof}

\begin{proof}[{\bf Proof of Theorem
1.15}]
The pinching condition in Theorem 1.15 implies that the equality holds in (50). So all inequalities leading to (49)
become equalities. Hence at every point, either $\mathring{Ric}$ is null, i.e., $M$ is Eninstein, or it has an eigenvalue of
multiplicity $(n-1)$ and another of multiplicity $1$. But (1) implies that $M$ is not Eninstein.  Since $M$ has harmonic curvature, and by the regularity result of DeTurck and Goldschmidt \cite{DG},
$M$ must be real analytic
in suitable (harmonic) local coordinates.

Suppose that the Ricci tensor has an eigenvalue of multiplicity $(n-1)$ and another of multiplicity $1$.
If the Ricci tensor is parallel, by the de Rham decomposition Theorem \cite{Dr},
$M$ is covered isometrically by the product of Einstein manifolds. We have $R=\sqrt{n(n-1)}|\mathring{Ric}|$. Since $M$ is conformally flat and has positive scalar curvature, then the only possibility is that $M$ is covered
isometrically by $\mathbb{S}^1\times \mathbb{S}^{n-1}$ with the product metric.

On the other hand, if the Ricci tensor is not parallel, by the  classification
result of Derdzi\'{n}ski (see Theorem 10 of \cite{D}, see also Theorem 3.2 of \cite{Ca}), this concludes the proof of Theorem 1.15.
\end{proof}

\end{document}